\documentclass[11pt]{article}
\usepackage[russian]{babel}
\usepackage[cp1251]{inputenc}
\usepackage{amsmath,latexsym}
\usepackage[psamsfonts]{amssymb}
\usepackage[dvips]{graphicx}
\usepackage[mathcal]{euscript}
\begin{document}
УДК. 517.954, 517.927.25
\begin{center}
\textbf{ ОБ УСЛОВИЯХ РАЗРЕШИМОСТИ  КРАЕВЫХ ЗАДАЧ\\ ДЛЯ ОДНОГО УРАВНЕНИЯ
ВЫСОКОГО ПОРЯДКА\\ С ПЕРЕМЕННЫМИ КОЭФФИЦИЕНТАМИ
 }\\
\textbf{Б.Ю.Иргашев}\\
Наманганский инженено-строительный институт.г.Наманган.Узбекистан.\\
Институт Математики Республики Узбекистан.\\
bahromirgasev@gmail.com
\end{center}

\textbf{Аннотация.} В работе для уравнения четного порядка с переменными коэффициентами изучена задача типа Дирихле. Дан критерий единственности решения. Решение построено в виде ряда Фурье . При обосновании сходимости ряда возникает проблема малых знаменателей. Получены достаточные условия отделимости знаменателя от нуля. Показана , что на разрешимость задачи влияет не только размерность прямоугольника, но и порядки задаваемых производных на нижней границе прямоугольника.

 \textbf{Ключевые слова.} Уравнение четного порядка, разрывный коэффициент, самосопряженная задача, собственное значение, собственная функция, определитель Вандермонда,малые знаменатели, единственность, ряд, равномерная сходимость, существование.

\textbf{1.Введение и постановка задачи.} В области $\Omega  = {\Omega _x} \times {\Omega _y},$ ${\Omega _x} = \left\{ {x:0 < x < \pi } \right\},$ ${\Omega _y} = \left\{ {y:\, - a < y < a} \right\}$ рассмотрим уравнение в частных производных
\[Lu \equiv l\left( {u\left( {x,y} \right)} \right) + {\left( { - 1} \right)^n}\left( {{\text{sgn}}y} \right){\mkern 1mu} D_y^{2n}u\left( {x,y} \right) = 0,\eqno(1)\]
где
\[l\left( {u\left( {x,y} \right)} \right) = {\left( { - 1} \right)^s}\frac{{{\partial ^{2s}}u\left( {x,y} \right)}}{{\partial {x^{2s}}}} + {\left( {{p_{s - 1}}\left( x \right)\frac{{{\partial ^{s - 1}}u\left( {x,y} \right)}}{{\partial {x^{s - 1}}}}} \right)^{\left( {s - 1} \right)}} + ...\]
\[ + {\left( {{p_1}\left( x \right)\frac{{\partial u\left( {x,y} \right)}}{{\partial x}}} \right)^\prime } + {p_0}\left( x \right)u\left( {x,y} \right),\]
\[{p_j}\left( x \right) \in {C^\infty }\left( {\overline {{\Omega _x}} } \right),\,\,j = 0,1,...,s - 1,\]
$$p_j^{\left( {2i + 1} \right)}\left( 0 \right) = p_j^{\left( {2i + 1} \right)}\left( \pi  \right) = 0,\,j = 0,1,...,s - 1,\,i \in N,$$
\[\frac{s}{n} = b \in N,\,D_t^ku = \frac{{{\partial ^k}u}}{{\partial {t^k}}},\,k \in N.\]
Пусть ${\Omega _ + } = \Omega  \cap \left( {y > 0} \right),{\mkern 1mu} {\mkern 1mu} {\Omega _ - } = \Omega  \cap \left( {y < 0} \right).$  Для уравнения  (1) рассмотрим следующую  задачу типа задачи Дирихле.

\textbf{Задача D.}  Найти функцию $u(x,y)$  с условиями
\[Lu\left( {x,y} \right) \equiv 0,{\mkern 1mu} {\mkern 1mu} {\mkern 1mu} {\mkern 1mu} \left( {x,y} \right) \in {\Omega _ + } \cup {\Omega _ - },\]
\[u\left( {x,y} \right) \in C_{x,y}^{2s - 1,2n - 1}\left( {\bar \Omega } \right) \cap C_{x,y}^{2s,2n}\left( {{\Omega _ + } \cup {\Omega _ - }} \right),\eqno(2)\]
\[D_x^{2j}u\left( {0,y} \right) = D_x^{2j}u\left( {\pi ,y} \right) = 0,{\mkern 1mu} {\mkern 1mu} \, - a \leqslant y \leqslant a,\eqno(3)\]
\[D_y^{q + \gamma r}u\left( {x, - a} \right) = {\varphi _r}\left( x \right),{\mkern 1mu} {\mkern 1mu} {\mkern 1mu} D_y^{\chi  + \delta r}u\left( {x,a} \right) = {\psi _r}\left( x \right),{\mkern 1mu} {\mkern 1mu} 0 \leqslant x \leqslant l,\eqno(4)\]
где
$${\varphi _r}\left( x \right),{\psi _r}\left( x \right) \in {C^{2s}}\left( {\overline {{\Omega _x}} } \right),\,\varphi _r^{\left( {2j} \right)}\left( 0 \right) = \varphi _r^{\left( {2j} \right)}\left( \pi  \right) = \psi _r^{\left( {2j} \right)}\left( 0 \right) = \psi _r^{\left( {2j} \right)}\left( \pi  \right) = 0,$$
$$j = 0,...,s - 1;{\mkern 1mu} r = 0,1,...,n - 1;\gamma ,\delta  \in \left\{ {1;2} \right\}.$$
Если $\gamma  = \delta  = 1,$ то  $q,\chi  \in \left\{ {0,1,...,n} \right\};$ если $\gamma  = \delta  = 2,$  то  $q,\chi  \in \left\{ {0,1} \right\}.$

Уравнение (1) при $s=n=1$ в области ${\Omega _ - }$ представляет собой известное уравнение струны. Задача Дирихле для уравнения струны изучалась во многих работах, например [1]–[8]. Из указанных работ следует, что если отношение сторон $a/\pi $ прямоугольника ${\Omega _ - }$, где ищется решение задачи Дирихле для уравнения струны, является рациональным числом, то однородная задача Дирихле имеет нетривиальные решения. В работе [9] показано, что иррациональность отношения $a/\pi $ является необходимым и достаточным условием единственности решения задачи Дирихле для модельного уравнения (1), при $s=n$, в области ${\Omega _ - }$.

В области $\Omega$ уравнение (1), при $s=n=1$ есть известное уравнение Лаврентьева-Бицадзе, для которого некорректность задачи Дирихле в смешанной области было показано А.В.Бицадзе [10] . После этого начались исследования по нахождению областей, для которых задача Дирихле была бы корректной, например работы [11]-[13]. В работах [14]-[19] задачу Дирихле для уравнений смешанного типа, в прямоугольной области, изучали методом разделения переменных. В работе [14] доказана корректность задачи Дирихле для уравнения Лаврентьева—Бицадзе в прямоугольной области при условии, когда отношение сторон прямоугольника в гиперболической части является натуральным числом.
	 В работе К.Б.Сабитова [17]  при обосновании сходимости построенного ряда Фурье была обнаружена проблема малых знаменателей.
Краевые задачи для уравнений высокого порядка с разрывными коэффициентами изучались в работах М.М.Смирнова, В.И.Жегалова, К.Б.Сабитова. В этих работах исследовались уравнение состоящие из произведений операторов Лаврентьева-Бицадзе.

В данной работе для уравнения (1) исследуется задача типа задачи Дирихле. Дан критерий однозначной разрешимости поставленной задачи.Решение строится в виде ряда Фурье по собственным функциям одномерной задачи. При обосновании сходимости ряда возникает проблема малых знаменателей. Получены достаточные условия отделимости малого знаменателя от нуля. Показана, что в отличии от уравнений второго порядка разрешимость поставленной задачи зависит не только от геометрии области, но и от порядков производных задаваемых на нижней границе прямоугольника.

\textbf{2. Единственность решения.} Учитывая (3), прежде рассмотрим следующую задачу на собственные значения:
\[\left\{ \begin{gathered}
  l\left( {X\left( x \right)} \right) = \lambda X\left( x \right), \hfill \\
  {X^{\left( {2j} \right)}}\left( 0 \right) = {X^{\left( {2j} \right)}}\left( \pi  \right) = 0,\,j = 0,1,...,s - 1. \hfill \\
\end{gathered}  \right.\eqno(5)\]
Опишем некоторые свойства собственных значений и собственных функций задачи (5). Задача (5) является самосопряженной, поэтому она имеет не более чем счетное число собственных значений ${\lambda _k},\,k = 1,2,...,$  и ортонормированную  систему собственных функций  ${X_k}\left( x \right)$ . Далее будем считать, что
$${\lambda _k} > 0,\,k = 1,2,...\,.\eqno(6)$$
условие (6) будем выполняться, если например
\[l\left( u \right) = {\left( { - 1} \right)^s}\frac{{{\partial ^{2s}}u\left( {x,y} \right)}}{{\partial {x^{2s}}}} + {p_0}\left( x \right)u\left( {x,y} \right),\,\,\,{p_0}\left( x \right) \geqslant 0.\eqno(7)\]
Из условия (6) следует, что у задачи (5) существует единственная функция Грина $G\left( {x,\xi } \right) = G\left( {\xi ,x} \right),$  с помощью которой задачу (5) можно свести к интегральному уравнению с симметричным ядром
\[{X_k}\left( x \right) = {\lambda _k}\int\limits_0^\pi  {G\left( {x,\xi } \right){X_k}\left( \xi  \right)d\xi } ,\]
отсюда учитывая теорему Мерсера (см. [24],стр.165) имеем
\[\sum\limits_{k = 1}^\infty  {\frac{{X_k^2\left( x \right)}}{{{\lambda _k}}}}  \leqslant G\left( {x,x} \right) < \infty .\eqno(8)\]
Известна (см. [23]) асимптотика собственных значений задачи (5):
\[{\lambda _k} = {k^{2s}} + {c_{ - 2s + 2}}{k^{2s - 2}} + ... + {c_0} + \frac{{{c_2}}}{{{k^2}}} + \frac{{{c_4}}}{{{k^4}}} + ...,\eqno(9)\]
где постоянные ${c_j}$  вычисляются определенным образом.

Пусть теперь $u(x,y)$ некоторое решение задачи $D.$ Рассмотрим его коэффициенты Фурье по системе собственных функций $X_{k}$
\[{Y_k}\left( y \right) = \int\limits_0^\pi  {u\left( {x,y} \right){X_k}\left( x \right)dx,} \]
введем следующую функцию
\[{Y_{k\varepsilon }}\left( y \right) = \int\limits_\varepsilon ^{\pi  - \varepsilon } {u\left( {x,y} \right){X_k}\left( x \right)dx,} \]
где $\epsilon >0$ -достаточно малое число, отсюда
\[Y_{k\varepsilon }^{\left( {2n} \right)}\left( y \right) = \int\limits_\varepsilon ^{\pi  - \varepsilon } {\frac{{{\partial ^{2n}}u\left( {x,y} \right)}}{{\partial {y^{2n}}}}{X_k}\left( x \right)dx}  = {\left( { - 1} \right)^{n + 1}}\left( {{\text{sgn}}y} \right)\int\limits_\varepsilon ^{\pi  - \varepsilon } {l\left( {u\left( {x,y} \right)} \right){X_k}\left( x \right)dx} .\]
Далее интегрируя по частям и переходя к пределу при $\varepsilon  \to  + 0$ , имеем
$$Y_k^{\left( {2n} \right)}\left( y \right) = {\left( { - 1} \right)^{n + 1}}\left( {{\text{sgn}}y} \right)\int\limits_0^\pi  {u\left( {x,y} \right)l\left( {{X_k}\left( x \right)} \right)dx}=$$
$$={\lambda _k}{\left( { - 1} \right)^{n + 1}}\left( {{\text{sgn}}y} \right)\int\limits_0^\pi  {u\left( {x,y} \right){X_k}\left( x \right)dx},$$
отсюда относительно переменной $y$   получим уравнение
\[D_y^{2n}{Y_k}\left( y \right){\text{ + sgn}}y{\left( { - 1} \right)^n}\lambda {Y_k}\left( y \right) = 0.\]
Рассмотрим отдельно случаи четного и нечетного $n$ . Пусть $n=2m$ (случай нечетного $n$  рассматривается аналогично) . Учитывая условия (2) и (4), получим следующую задачу:
\[\left\{ \begin{gathered}
  Y_k^{\left( {4m} \right)}\left( y \right) + \operatorname{sgn} y\lambda {Y_k}\left( y \right) = 0, \hfill \\
  Y_k^{\left( {q + \gamma j} \right)}\left( { - a} \right) = {\varphi _{jk}}, \hfill \\
  Y_k^{\left( {\chi  + \delta j} \right)}\left( a \right) = {\psi _{jk}}, \hfill \\
  Y_k^{\left( r \right)}\left( { + 0} \right) = Y_k^{\left( r \right)}\left( { - 0} \right),{\mkern 1mu} {\kern 1pt} j = \overline {0,2m - 1} ,{\mkern 1mu} {\kern 1pt} {\mkern 1mu} {\kern 1pt} r = \overline {0,(4m - 1)} . \hfill \\
\end{gathered}  \right.\eqno(10)\]
где
\[{\varphi _{jk}} = \int\limits_0^l {{\varphi _j}\left( x \right){X_k}\left( x \right)dx,{\mkern 1mu} {\kern 1pt} {\mkern 1mu} {\kern 1pt} } {\psi _{jk}} = \int\limits_0^l {{\psi _j}\left( x \right){X_k}\left( x \right)dx} .\]
Отметим, что возможность разложения граничных функций ${\varphi _j}\left( x \right),{\psi _j}\left( x \right)$ в ряд Фурье по системе собственных функций ${X_k}\left( x \right)$, следует из условий наложенных на эти функции и теоремы Гильберта-Шмидта.

При $y > 0$ общее решение уравнения (10) имеет вид
\[{Y_k}\left( y \right) = \sum\limits_{r = 0}^{2m - 1} {\left( {c_r^1{Y_{1r}}(y) + c_r^2{Y_{2r}}\left( y \right)} \right)} ,\]
здесь
\[{Y_{1r}}\left( y \right) = {e^{{\alpha _r}y}}\cos {\beta _r}y,{\mkern 1mu} {\kern 1pt} {\mkern 1mu} {\kern 1pt} {Y_{2r}}\left( y \right) = {e^{{\alpha _r}y}}\sin {\beta _r}y,\]
\[{\alpha _r} = \sqrt[{4m}]{{{\lambda _k}}}\cos {\theta _r},{\beta _r} = \sqrt[{4m}]{{{\lambda _k}}}\sin {\theta _r},{\mkern 1mu} {\kern 1pt} {\theta _r} = \frac{\pi }{{4m}}\left( {1 + 2r} \right),\]
\[{\left( {{Y_{1r}}\left( y \right)} \right)^{\left( j \right)}} = \sqrt[{4m}]{{\lambda _k^j}}{e^{{\alpha _r}y}}\cos \left( {{\beta _r}y + j{\theta _r}} \right),\]
\[{\left( {{Y_{2r}}\left( y \right)} \right)^{\left( j \right)}} = \sqrt[{4m}]{{\lambda _k^j}}{e^{{\alpha _r}y}}sin\left( {{\beta _r}y + j{\theta _r}} \right),\]
\[r = \overline {0,(2m - 1)} ,{\mkern 1mu} {\kern 1pt} j = \overline {0,4m - 1} ,{\mkern 1mu} {\kern 1pt} {\mkern 1mu} {\kern 1pt} {\alpha _r} > 0,{\mkern 1mu} {\kern 1pt} r = \overline {0,(m - 1)} .\]
При $y < 0$ имеем
\[{Y_k}\left( y \right) = {d_0}{e^{\sqrt[{4m}]{{{\lambda _k}}}y}} + \sum\limits_{r = 1}^{2m - 1} {{e^{{\mu _r}y}}\left( {d_r^1\cos {\nu _r}y + d_r^2\sin {\nu _r}y} \right)}  + {d_{2m}}{e^{ - \sqrt[{4m}]{{{\lambda _k}}}y}},\]
\[Y_k^{\left( {q + \gamma j} \right)}\left( y \right) = \sqrt[{4m}]{{\lambda _k^{q + \gamma j}}}\left( {{d_0}{e^{\sqrt[{4m}]{{{\lambda _k}}}y}} + } \right.\sum\limits_{r = 1}^{2m - 1} {{e^{{\mu _r}y}}d_r^1\cos \left( {{\nu _r}y + \left( {q + \gamma j} \right){\sigma _r}} \right) + } \]
\[ + \left. {\sum\limits_{r = 1}^{2m - 1} {{e^{{\mu _r}y}}d_r^2\sin \left( {{\nu _r}y + \left( {q + \gamma j} \right){\sigma _r}} \right)}  + {{\left( { - 1} \right)}^{\left( {q + \gamma j} \right)}}{d_{2m}}{e^{ - \sqrt[{4m}]{{{\lambda _k}}}y}}} \right),\]
где
\[{\mu _r} = \sqrt[{4m}]{{{\lambda _k}}}\cos {\sigma _r},{\nu _r} = \sqrt[{4m}]{{{\lambda _k}}}\sin {\sigma _r},{\mkern 1mu} {\kern 1pt} {\sigma _r} = \frac{{\pi r}}{{2m}},\]
\[r = \overline {0,(2m - 1)} ,{\mu _r} < 0,r = \overline {m + 1,2m} ,{\mkern 1mu} {\kern 1pt} {\mu _m} = 0.\]
Удовлетворив краевым условиям задачи (10) получим систему алгебраических уравнений
\[\left\{ \begin{gathered}
  \sum\limits_{r = 0}^{2m - 1} {{e^{{\alpha _r}y}}\left( {c_r^1\cos \left( {{\beta _r}a + \left( {\chi  + \delta j} \right){\theta _r}} \right) + } \right.}  \hfill \\
   + c_r^2sin\left( {{\beta _r}a + \left( {\chi  + \delta j} \right){\theta _r}} \right) = {\left( {\frac{1}{{\sqrt[{4m}]{{{\lambda _k}}}}}} \right)^{\chi  + \delta j}}{\psi _{jk}}, \hfill \\
  \sum\limits_{r = 0}^{2m - 1} {{e^{ - a{\mu _r}}}\left( {d_r^1\cos \left( { - {\nu _r}a + \left( {q + \gamma j} \right){\sigma _r})} \right)} \right. + }  \hfill \\
   + \left. {d_r^2\sin \left( { - {\nu _r}a + \left( {q + \gamma j} \right){\sigma _r}} \right)} \right) = {\left( {\frac{1}{{\sqrt[{4m}]{{{\lambda _k}}}}}} \right)^{q + j\gamma }}{\varphi _{jk}}, \hfill \\
  \sum\limits_{r = 0}^{2m - 1} {\left( {c_r^1\cos \left( {l{\theta _r}} \right) + c_r^2\sin \left( {l{\theta _r}} \right)} \right)}  =  \hfill \\
   = {d_0} + \sum\limits_{r = 1}^{2m - 1} {\left( {d_r^1\cos \left( {l{\sigma _r}} \right) + d_r^2\sin \left( {l{\sigma _r}} \right)} \right) + {{\left( { - 1} \right)}^l}{d_{2m}},}  \hfill \\
  j = \overline {0,2m - 1} ,{\mkern 1mu} {\kern 1pt} {\mkern 1mu} {\kern 1pt} {\mkern 1mu} {\kern 1pt} l = \overline {0,(4m - 1)} . \hfill \\
\end{gathered}  \right.\eqno(11)\]
Введем следующие обозначения :
\[{\omega _{j,r}} = {\beta _r}a + \left( {\chi  + \delta j} \right){\theta _r},{\mkern 1mu} {\kern 1pt} {\tau _{j,r}} =  - {\nu _r}a + \left( {q + \gamma j} \right){\sigma _r},\]
$$\alpha  = {\alpha _0} + {\alpha _1} + ... + {\alpha _{m - 1}}, \mu  = {\mu _{m + 1}} + ... + {\mu _{2m - 1}},$$
\[A_{2m}^ +  = \left( {\begin{array}{*{20}{c}}
  {{e^{{\alpha _0}a}}\cos {\omega _{0,0}}}&{{e^{{\alpha _0}a}}sin{\omega _{0,0}}}&.&.&{{e^{{\alpha _{m - 1}}a}}\sin {\omega _{0,m - 1}}} \\
  {{e^{{\alpha _0}a}}\cos {\omega _{1,0}}}&{{e^{{\alpha _0}a}}sin{\omega _{1,0}}}&.&.&{{e^{{\alpha _{m - 1}}a}}\sin {\omega _{1,m - 1}}} \\
  .&.&.&.&. \\
  {{e^{{\alpha _0}a}}\cos {\omega _{2m - 1,0}}}&{{e^{{\alpha _0}a}}sin{\omega _{2m - 1,0}}}&.&.&{{e^{{\alpha _{m - 1}}a}}\sin {\omega _{2m - 1,m - 1}}}
\end{array}} \right),\]
\[A_{2m,2m}^ -  = \left( {\begin{array}{*{20}{c}}
  {{e^{{\alpha _m}a}}\cos {\omega _{0,m}}}&{{e^{{\alpha _m}a}}sin{\omega _{0,m}}}&.&.&{{e^{{\alpha _{2m - 1}}a}}\sin {\omega _{0,2m - 1}}} \\
  {{e^{{\alpha _m}a}}\cos {\omega _{1,m}}}&{{e^{{\alpha _m}a}}sin{\omega _{1,m}}}&.&.&{{e^{{\alpha _{2m - 1}}a}}\sin {\omega _{1,2m - 1}}} \\
  .&.&.&.&. \\
  {{e^{{\alpha _m}a}}\cos {\omega _{2m - 1,m}}}&{{e^{{\alpha _m}a}}sin{\omega _{2m - 1,m}}}&.&.&{{e^{{\alpha _{2m - 1}}a}}\sin {\omega _{2m - 1,2m - 1}}}
\end{array}} \right),\]
\[B_{2m,2m - 1}^ +  = \left( {\begin{array}{*{20}{c}}
  {{e^{ - {\mu _{m + 1}}a}}\cos {\tau _{0,m + 1}}}&.&{{{\left( { - 1} \right)}^q}{e^{\sqrt[{4m}]{{{\lambda _k}}}a}}} \\
  .&.&. \\
  {{e^{ - {\mu _{m + 1}}a}}\cos {\tau _{2m - 1,m + 1}}}&.&{{{\left( { - 1} \right)}^{q + \gamma \left( {2m - 1} \right)}}{e^{\sqrt[{4m}]{{{\lambda _k}}}a}}}
\end{array}} \right),\]
или  в компактной записи
\[B_{2m,2m - 1}^ +  = \left( {{e^{ - {\mu _r}a}}\cos {\tau _{j,r}},{e^{ - {\mu _r}a}}\sin {\tau _{j,r}},{{\left( { - 1} \right)}^{q + \gamma j}}{e^{\sqrt[{4m}]{{{\lambda _k}}}a}}} \right)_{r = \overline {m + 1,2m - 1} }^{j = \overline {0,2m - 1} },\]
здесь верхний индекс означает строку , а нижний столбец. Аналогично
\[B_{2m,2m - 1}^ -  = \left( {{e^{ - \sqrt[{4m}]{{{\lambda _k}}}a}},{e^{ - {\mu _r}a}}\cos {\tau _{j,r}},{e^{ - {\mu _r}a}}\sin {\tau _{j,r}}} \right)_{r = \overline {1,m - 1} }^{j = \overline {0,2m - 1} },\]
\[C_{4m,2m}^ +  = \left( {\cos j{\theta _r},\sin j{\theta _r}} \right)_{r = \overline {0,m - 1} }^{j = \overline {0,4m - 1} },{\mkern 1mu} {\kern 1pt} C_{4m,2m}^ -  = \left( {\cos j{\theta _r},\sin j{\theta _r}} \right)_{r = \overline {m,2m - 1} }^{j = \overline {0,4m - 1} },\]
\[D_{4m,2m - 1}^ +  = \left( { - \cos j{\sigma _r}, - \sin j{\sigma _r},{{\left( { - 1} \right)}^j}} \right)_{r = \overline {m + 1,2m - 1} }^{j = \overline {0,4m - 1} },\]
\[D_{4m,2m - 1}^ -  = \left( {1, - \cos j{\sigma _r}, - \sin j{\sigma _r}} \right)_{r = \overline {1,m - 1} }^{j = \overline {0,4m - 1} }.\]
Далее формально имеем
\[B_{2m,2}^0 = \left( {\begin{array}{*{20}{c}}
  {\cos {\tau _{0,m}}}&{\sin {\tau _{0,m}}} \\
  {\cos {\tau _{1,m}}}&{\sin {\tau _{1,m}}} \\
  .&. \\
  {\cos {\tau _{2m - 1,m}}}&{\sin {\tau _{2m - 1,m}}}
\end{array}} \right) = \frac{i}{2}\left( {\begin{array}{*{20}{c}}
  {{e^{i{\tau _{0,m}}}}}&{{e^{ - i{\tau _{0,m}}}}} \\
  .&. \\
  {{e^{i{\tau _{2m - 1,m}}}}}&{{e^{ - i{\tau _{2m - 1,m}}}}}
\end{array}} \right) = \frac{i}{2}\left( {\begin{array}{*{20}{c}}
  {B_{2m,1}^ + }&{B_{2m,1}^ - }
\end{array}} \right),\]
\[D_{4m,2}^0 = \left( {\begin{array}{*{20}{c}}
  {\cos 0}&{\sin 0} \\
  {\cos \frac{\pi }{2}}&{\sin \frac{\pi }{2}} \\
  .&. \\
  {\cos \left( {4m - 1} \right)\frac{\pi }{2}}&{\sin \left( {4m - 1} \right)\frac{\pi }{2}}
\end{array}} \right) = \]
\[ = \frac{i}{2}\left( {\begin{array}{*{20}{c}}
  {{e^{i0 \cdot \frac{\pi }{2}}}}&{{e^{ - i0 \cdot \frac{\pi }{2}}}} \\
  .&. \\
  {{e^{i\left( {4m - 1} \right)\frac{\pi }{2}}}}&{{e^{ - i\left( {4m - 1} \right)\frac{\pi }{2}}}}
\end{array}} \right) = \frac{i}{2}\left( {\begin{array}{*{20}{c}}
  {D_{4m,1}^ + }&{D_{4m,1}^ - }
\end{array}} \right).\]
В этих обозначениях основной определитель системы (11), ${\Delta _1}\left( k \right)$ запишется в виде
\[{\Delta _1}\left( k \right) = \frac{i}{2}\det \left( {\begin{array}{*{20}{c}}
  {A_{2m,2m}^ + }&{A_{2m,2m}^ - }&0&0&0&0 \\
  0&0&{B_{2m,2m - 1}^ - }&{B_{2m,2m - 1}^ + }&{B_{2m,1}^ + }&{B_{2m,1}^ - } \\
  {C_{4m,2m}^ + }&{C_{4m,2m}^ - }&{D_{4m,2m - 1}^ - }&{D_{4m,2m - 1}^ + }&{D_{4m,1}^ + }&{D_{4m,1}^ - }
\end{array}} \right).\]
Найдем асимптотику определителя ${\Delta _1}\left( k \right)$  , при больших значениях $k$. Для этого вычислим слагаемое, куда входит экспонента с наибольшей положительной степенью. С точностью до знака он имеет вид
\[{\Delta _2}\left( k \right) = \frac{i}{2}\left| {A_{2m,2m}^ + } \right|\left( {\left| {B_{2m,2m - 1}^ + B_{2m,1}^ + } \right| \cdot \left| {C_{4m,2m}^ - D_{4m,2m - 1}^ - D_{4m,1}^ - } \right| - } \right.\]
\[ - \left. {\left| {B_{2m,2m - 1}^ + B_{2m,1}^ - } \right| \cdot \left| {C_{4m,2m}^ - D_{4m,2m - 1}^ - D_{4m,1}^ + } \right|} \right).\]
Перейдем к вычислениям. Действовать будем следующим образом, используя формулу Эйлера ${e^{it}} = \cos t + i\sin t$, приведем вычисляемые определители к определителя Вандермонда
\[\left| {\begin{array}{*{20}{c}}
  1&1&.&.&1 \\
  {{x_1}}&{{x_2}}&.&.&{{x_n}} \\
  .&.&.&.&. \\
  {x_1^{n - 1}}&{x_2^{n - 1}}&.&.&{x_n^{n - 1}}
\end{array}} \right| = \prod\limits_{j > i} {\left( {{x_j} - {x_i}} \right)} .\]
Имеем
\[\left| {A_{2m,2m}^ + } \right| = {\left( {\frac{i}{2}} \right)^m}{e^{2\alpha a}}\left| {\begin{array}{*{20}{c}}
  {{e^{i{\omega _{0,0}}}}}&{{e^{ - i{\omega _{0,0}}}}}&.&{{e^{i{\omega _{0,m - 1}}}}}&{{e^{ - i{\omega _{0,m - 1}}}}} \\
  {{e^{i{\omega _{1,0}}}}}&{{e^{ - i{\omega _{1,0}}}}}&.&{{e^{i{\omega _{1,m - 1}}}}}&{{e^{ - i{\omega _{1,m - 1}}}}} \\
  .&.&.&.&. \\
  {{e^{i{\omega _{2m - 1,0}}}}}&{{e^{ - i{\omega _{2m - 1,0}}}}}&.&{{e^{i{\omega _{2m - 1,m - 1}}}}}&{{e^{ - i{\omega _{2m - 1,m - 1}}}}}
\end{array}} \right| = \]
\[ = {\left( {\frac{i}{2}} \right)^m}{e^{2\alpha a}}\left| {\begin{array}{*{20}{c}}
  1&1&.&1&1 \\
  {{e^{i\delta {\theta _0}}}}&{{e^{ - i\delta {\theta _0}}}}&.&{{e^{i\delta {\theta _{m - 1}}}}}&{{e^{ - i\delta {\theta _{m - 1}}}}} \\
  .&.&.&.&. \\
  {{e^{i\delta \left( {2m - 1} \right){\theta _0}}}}&{{e^{ - i\delta \left( {2m - 1} \right){\theta _0}}}}&.&{{e^{i\delta \left( {2m - 1} \right){\theta _{m - 1}}}}}&{{e^{ - i\delta \left( {2m - 1} \right){\theta _{m - 1}}}}}
\end{array}} \right| = \]
\[ = {\left( {\frac{i}{2}} \right)^m}{e^{2\alpha a}}\prod\limits_{j = 0}^{m - 1} {\left( { - 2i\sin \delta {\theta _j}} \right)}  \cdot \]
\[ \cdot \prod\limits_{0 = s < j = m - 1}^{} {4\left( {1 - \cos \delta \left( {{\theta _j} - {\theta _s}} \right)} \right)\left( {1 - \cos \delta \left( {{\theta _j} + {\theta _s}} \right)} \right)}  \ne 0,\]
Действую аналогичным образом придем к следующим результатам:
\[\det \left( {B_{2m,2m - 1}^ + ,B_{2m,1}^ + } \right) = \]
\[ = {e^{ - 2a\mu  + \sqrt[{4m}]{{{\lambda _k}}}a}}{e^{ - \sqrt[{4m}]{{{\lambda _k}}}ai}}{\left( {\frac{i}{2}} \right)^{m - 1}}{\left( { - i} \right)^q}{M_1}\left( {{i^\gamma } - {{\left( { - 1} \right)}^\gamma }} \right) \cdot \]
\[ \cdot \prod\limits_{r = m + 1}^{2m - 1} {\left( {{{\left( { - 1} \right)}^\gamma } + 1 - 2{i^\gamma }\cos \gamma {\sigma _r}} \right)} ,\]
здесь
\[{M_1} = \det \left( {{e^{i\gamma j{\sigma _r}}},{e^{ - i\gamma j{\sigma _r}}},{{\left( { - 1} \right)}^{\gamma j}}} \right)_{r = \overline {m + 1,2m - 1} }^{j = \overline {0,2m - 2} } \ne 0.\]
Далее
\[\det \left( {B_{2m,2m - 1}^ + ,B_{2m,1}^ - } \right) = \]
\[ = {e^{ - 2\mu a + \sqrt[{4m}]{{{\lambda _k}}}a}}{e^{i\sqrt[{4m}]{{{\lambda _k}}}a}}{\left( {\frac{i}{2}} \right)^{m - 1}}{i^q} \cdot \]
\[ \cdot {M_1}\left( {{{\left( { - i} \right)}^\gamma } - {{\left( { - 1} \right)}^\gamma }} \right)\prod\limits_{r = m + 1}^{2m - 1} {\left( {{{\left( { - 1} \right)}^\gamma } + 1 - 2{{\left( { - i} \right)}^\gamma }\cos \gamma {\sigma _r}} \right)} ,\]
\[\det \left( {C_{4m,2m}^ - D_{4m,2m - 1}^ - D_{4m,1}^ - } \right) = \]
\[ = {M_2}\left( {i + 1} \right)\prod\limits_{l = 1}^{m - 1} {\cos {\sigma _l}} \prod\limits_{r = m}^{2m - 1} {\cos {\theta _r}}  = {M_3}\left( {i + 1} \right),\]
где
\[{M_2} = \det \left( {{e^{ij{\theta _r}}},{e^{ - ij{\theta _r}}},1,{e^{ij{\sigma _l}}},{e^{ - ij{\sigma _l}}}} \right)_{r = \overline {m,2m - 1} ;l = \overline {1,m - 1} }^{j = \overline {0,4m - 2} } \ne 0.\]
\[{M_3} = {M_2}\prod\limits_{l = 1}^{m - 1} {\cos {\sigma _l}} \prod\limits_{r = m}^{2m - 1} {\cos {\theta _r}}  \ne 0.\]
Аналогично
\[\det \left( {C_{4m,2m}^ - D_{4m,2m - 1}^ - D_{4m,1}^ + } \right) = {M_3}\left( {i - 1} \right).\]
Далее все постоянные не зависящие от $k$  , будем обозначать одной буквой $L$. Учитывая это имеем
\[{\Delta _2}\left( k \right) = L{e^{2a\left( {\alpha  - \mu } \right) + \sqrt[{4m}]{{{\lambda _k}}}a}}{\Delta _3}\left( k \right),\]
здесь
\[{\Delta _3}\left( k \right) = \left( {{{\left( { - i} \right)}^q}\left( {{i^\gamma } - {{\left( { - 1} \right)}^\gamma }} \right)\prod\limits_{r = m + 1}^{2m - 1} {\left( {{{\left( { - 1} \right)}^\gamma } + 1 - 2{i^\gamma }\cos \gamma {\sigma _r}} \right)} \left( {i + 1} \right){e^{ - \sqrt[{4m}]{{{\lambda _k}}}ai}} - } \right.\]
\[\left. { - {i^q}\left( {{{\left( { - i} \right)}^\gamma } - {{\left( { - 1} \right)}^\gamma }} \right)\prod\limits_{r = m + 1}^{2m - 1} {\left( {{{\left( { - 1} \right)}^\gamma } + 1 - 2{{\left( { - i} \right)}^\gamma }\cos \gamma {\sigma _r}} \right)} \left( {i - 1} \right){e^{i\sqrt[{4m}]{{{\lambda _k}}}a}}} \right).\]
Рассмотрев все частные случаи приходим к следующему результату. Основной определитель системы (11) имеет вид
\[{\Delta _1}\left( k \right) = L{e^{2a\left( {\alpha  - \mu } \right) + \sqrt[{2n}]{{{\lambda _k}}}a}}\left( {{\Delta _4}\left( k \right) + {\Delta _5}\left( k \right)} \right),\]
где
\[{\Delta _4}\left( k \right) = \left\{ \begin{gathered}
  \sin \left( {\sqrt[{2n}]{{{\lambda _k}}}a + \frac{\pi }{2}} \right),\left\{ {2n = 8l + 4,{\mkern 1mu} {\kern 1pt} \gamma  = 1,{\mkern 1mu} {\kern 1pt} q = 2j + 1} \right\} \cup \left\{ {2n = 8l,{\mkern 1mu} {\kern 1pt} \gamma  = 1,{\mkern 1mu} {\kern 1pt} q = 2j} \right\}; \hfill \\
  \sin \sqrt[{2n}]{{{\lambda _k}}}a,\,\left\{ {2n = 8l + 4,{\mkern 1mu} {\kern 1pt} \gamma  = 1,{\mkern 1mu} {\kern 1pt} q = 2j} \right\} \cup \left\{ {2n = 8l,{\mkern 1mu} {\kern 1pt} \gamma  = 1,{\mkern 1mu} {\kern 1pt} q = 2j + 1} \right\}; \hfill \\
  \sin \left( {\sqrt[{2n}]{{{\lambda _k}}}a + \frac{\pi }{4}} \right),\left\{ {2n = 4l,{\mkern 1mu} {\kern 1pt} \gamma  = 2,q = 2j} \right\}; \hfill \\
  \sin \left( {\sqrt[{2n}]{{{\lambda _k}}}a + \frac{{3\pi }}{4}} \right),\left\{ {2n = 4l,{\mkern 1mu} {\kern 1pt} \gamma  = 2,q = 2j + 1} \right\}; \hfill \\
  l,j \in N \cup \left\{ 0 \right\}, \hfill \\
\end{gathered}  \right.\]
\[\mathop {\lim }\limits_{k \to \infty } {\Delta _5}\left( k \right) = 0.\]
В случаи когда $n=2m+1$ получим следующий результат:
\[{\Delta _1}\left( k \right) = L{e^{2a\left( {\alpha  - \mu } \right) + \sqrt[{2n}]{{{\lambda _k}}}a}}\left( {{\Delta _4}\left( k \right) + {\Delta _5}\left( k \right)} \right),\]
где
\[{\Delta _4}\left( k \right) = \left\{ \begin{gathered}
  \sin \left( {\sqrt[{2n}]{{{\lambda _k}}}a + \frac{\pi }{4}} \right),\left\{ {2n = 8l + 2,\gamma  = 1,q = 2j} \right\}, \hfill \\
  \left\{ {2n = 8l + 6,\gamma  = 1,q = 2j + 1} \right\},\left\{ {2n = 4l + 2,\gamma  = 2,q = 2j} \right\}; \hfill \\
  \sin \left( {\sqrt[{2n}]{{{\lambda _k}}}a + \frac{{3\pi }}{4}} \right),\left\{ {2n = 8l + 6,\gamma  = 1,q = 2j} \right\}, \hfill \\
  \left\{ {2n = 8l + 2,\gamma  = 1,q = 2j + 1} \right\},\left\{ {2n = 4l + 2,\gamma  = 2,q = 2j + 1} \right\}, \hfill \\
  l,j \in N \cup \left\{ 0 \right\}, \hfill \\
\end{gathered}  \right.\]
\[\mathop {\lim }\limits_{k \to \infty } {\Delta _5}\left( k \right) = 0.\]
Справедлива следующая теорема единственности.

\textbf{Теорема 1.} Если решение задачи $D$ существует , то оно единственно тогда и только тогда, когда определитель системы (11) отличен от нуля для всех значений $k$.

Доказательство следует из полноты системы собственных функций задачи (5) (см.[25]).

\textbf{3. Существование решения.}

Сначала получим некоторые оценки для ${Y_k}\left( y \right)$ . Справедлива лемма.

\textbf{Лемма 1.} Для функции ${Y_k}\left( y \right)$ и её производных, при достаточно больших значениях $k$ , справедливы оценки

\[\left| {Y_k^{(j)}\left( y \right)} \right| \leqslant L\frac{{\sqrt[{2n}]{{\lambda _k^j}}\sum\limits_{r = 0}^{n - 1} {\left\{ {\left| {{\varphi _{rk}}} \right| + \left| {{\psi _{rk}}} \right|} \right\}} }}{{\left| {{\Delta _4}\left( k \right) + {\Delta _5}\left( k \right)} \right|}},j = 0,1,...,2n.\]

\textbf{Доказательство.} Пусть $ n = 2m,y > 0 $  (другие случаи рассматриваются аналогично),  нетрудно показать, что
\[\left| {Y_k^{\left( j \right)}\left( y \right)} \right| \leqslant L\sqrt[{2n}]{{\lambda _k^j}}\sum\limits_{r = 0}^{2m - 1} {{e^{{\alpha _r}y}}\left( {\left| {c_r^1} \right| + \left| {c_r^2} \right|} \right)} ,j = 0,1,...,4m,\]
поэтому достаточно доказать оценку для $j = 0.$  Имеем
\[{e^{{\alpha _0}y}}\left| {c_0^1} \right| = {e^{{\alpha _0}y}}\left| {\frac{{{\Delta _0}\left( k \right)}}{{{\Delta _1}\left( k \right)}}} \right| \leqslant {e^{{\alpha _0}a}}\frac{{\sum\limits_{r = 0}^{2m - 1} {\left\{ {\left| {{\varphi _{rk}}} \right| + \left| {{\psi _{rk}}} \right|} \right\}}  \cdot O\left( {{e^{2a\left( {\alpha  - \mu } \right) + \sqrt[{4m}]{{{\lambda _k}}}a - {\alpha _0}a}}} \right)}}{{{e^{2a\left( {\alpha  - \mu } \right) + \sqrt[{4m}]{{{\lambda _k}}}a}}\left| {{\Delta _4}\left( k \right) + {\Delta _5}\left( k \right)} \right|}} \leqslant \]
\[ \leqslant L\frac{{\sum\limits_{r = 0}^{2m - 1} {\left\{ {\left| {{\varphi _{rk}}} \right| + \left| {{\psi _{rk}}} \right|} \right\}} }}{{\left| {{\Delta _4}\left( k \right) + {\Delta _5}\left( k \right)} \right|}},\]
здесь ${\Delta _0}\left( k \right)$ - определитель матрицы, полученный заменой первого столбца основной матрицы системы (11), правой частью системы (11). Справедливость полученной оценки для других слагаемых показывается аналогично.	\textbf{Лемма 1 доказана}.

Нужно теперь найти условия, при которых выражение  $\left| {{\Delta _4}\left( k \right) + {\Delta _5}\left( k \right)} \right|$  отделяется от нуля, т.е. начиная с некоторого номера $k$ должно выполняться оценка
\[\left| {{\Delta _4}\left( k \right) + {\Delta _5}\left( k \right)} \right| > {\delta _1} > 0.\]
Т.к. $\mathop {\lim }\limits_{k \to \infty } {\Delta _5}\left( k \right) = 0,$  то достаточно найти условия, при которых $\left| {{\Delta _4}\left( k \right)} \right| > {\delta _1} > 0,{\mkern 1mu} {\kern 1pt} k \to \infty .$  Справедлива лемма.

\textbf{Лемма 2.} Пусть собственные значения задачи (5) имеют асимптотику
\[{\lambda _k} = {k^{2s}} + O\left( {{k^r}} \right),\,r < 2s - b,\eqno(12)\]
тогда для справедливости оценки
\[\left| {{\Delta _4}\left( k \right)} \right| \geqslant {\delta _1} > 0,{\mkern 1mu} {\kern 1pt} k \to  + \infty ,\eqno(13)\]
достаточно выполнения одного из двух условий:\\
1. ${\Delta _4}\left( k \right) = \sin \left( {\sqrt[{2n}]{{{\lambda _k}}}a + \frac{\pi }{4}} \right)$  или  ${\Delta _4}\left( k \right) = \sin \left( {\sqrt[{2n}]{{{\lambda _k}}}a + \frac{{3\pi }}{4}} \right);$ $\frac{a}{\pi } \in N,$   либо   $\frac{a}{\pi } = \frac{s}{t}{\mkern 1mu} {\kern 1pt} {\mkern 1mu} {\kern 1pt} {\mkern 1mu} {\kern 1pt} \left( {\frac{a}{\pi } \notin N} \right),$ $s,t \in N,{\mkern 1mu} {\kern 1pt} {\mkern 1mu} {\kern 1pt} \left( {s,t} \right) = 1,{\mkern 1mu} $ $t$ не кратно 4.\\
2. ${\Delta _4}\left( k \right) = \sin \left( {\sqrt[{2n}]{{{\lambda _k}}}a + \frac{\pi }{2}} \right);$ $\frac{a}{\pi } \in N,$ либо
$\frac{a}{\pi } = \frac{s}{t}{\mkern 1mu} {\kern 1pt} {\mkern 1mu} {\kern 1pt} {\mkern 1mu} {\kern 1pt} \left( {\frac{a}{\pi } \notin N} \right),$ $s,t \in N,{\mkern 1mu} {\kern 1pt} {\mkern 1mu} {\kern 1pt} \left( {s,t} \right) = 1,{\mkern 1mu} {\kern 1pt} {\mkern 1mu} {\kern 1pt} {\mkern 1mu} {\kern 1pt} \left( {t,2} \right) = 1.$

\textbf{Доказательство.} Заметим, что при выполнении условия (12)  имеем
\[\sqrt[{2n}]{{{\lambda _k}}} = {\left( {{k^{2s}} + O\left( {{k^r}} \right)} \right)^{\frac{1}{{2n}}}} = {k^{\frac{s}{n}}}{\left( {1 + O\left( {{k^{r - 2s}}} \right)} \right)^{\frac{1}{{2n}}}} = {k^b} + O\left( {{k^{r + b - 2s}}} \right).\]
Перейдем к доказательству оценки (13). Пусть  в условии 2 ( 1. проверяется аналогично),  $\frac{a}{\pi } \in N$ , тогда
\[\left| {\sin \left( {\sqrt[{2n}]{{{\lambda _k}}}\pi \frac{a}{\pi } + \frac{\pi }{2}} \right)} \right| = \left| {\cos \left( {\sqrt[{2n}]{{{\lambda _k}}}\pi \frac{a}{\pi }} \right)} \right| = \]
\[ = \left| {\cos \left( {{k^b}\frac{a}{\pi }\pi  + aO\left( {{k^{r + b - 2s}}} \right)} \right)} \right| = \left| {\cos \left( {aO\left( {{k^{r + b - 2s}}} \right)} \right)} \right| > {\delta _1} > 0,k \to  + \infty ,\]
т.е. (13) выполнено.

Пусть теперь $\frac{a}{\pi } = \frac{s}{t}{\mkern 1mu} {\kern 1pt} {\mkern 1mu} {\kern 1pt} {\mkern 1mu} {\kern 1pt} \left( {\frac{a}{\pi } \notin N} \right),$ тогда
\[\left| {\sin \left( {\sqrt[{2n}]{{{\lambda _k}}}\pi \frac{a}{\pi } + \frac{\pi }{2}} \right)} \right| = \left| {\sin \left( {\pi \frac{{s{k^b}}}{t} + \frac{\pi }{2} + aO\left( {{k^{r + b - 2s}}} \right)} \right)} \right| \geqslant \]
\[ \geqslant \left| {\sin \left( {\pi \frac{{s{k^b}}}{t} + \frac{\pi }{2}} \right)} \right|\left| {\cos \left( {aO\left( {{k^{r + b - 2s}}} \right)} \right)} \right| - \left| {\sin \left( {aO\left( {{k^{r + b - 2s}}} \right)} \right)} \right|\left| {sin\left( {\pi \frac{{s{k^b}}}{t}} \right)} \right|.\]
Осталось показать, что
\[\left| {\sin \left( {\pi \frac{{s{k^b}}}{t} + \frac{\pi }{2}} \right)} \right| \geqslant {\delta _1} > 0,\]
действительно, имеем  $\frac{{{k^b}a}}{\pi } = \frac{{{k^b}s}}{t} = {k_1} + \frac{{{k_2}}}{t},$  где  ${k_1},{k_2} \in N,{\mkern 1mu} {\kern 1pt} {\mkern 1mu} {\kern 1pt} 1 \leqslant {k_2} \leqslant t - 1.$ Значит
\[\left| {\sin \left( {\pi \frac{{s{k^b}}}{t} + \frac{\pi }{2}} \right)} \right| = \left| {\sin \pi \left( {\frac{{{k_2}}}{t} + \frac{1}{2}} \right)} \right|.\]
Т.к.
\[t \ne 2{k_2} \Rightarrow \frac{{{k_2}}}{t} \ne \frac{1}{2} \Rightarrow \frac{{{k_2}}}{t} + \frac{1}{2} \ne 1 \Rightarrow \left| {\sin \pi \left( {\frac{{{k_2}}}{t} + \frac{1}{2}} \right)} \right| > 0 \Rightarrow \]
\[\left| {\sin \left( {\pi \frac{{s{k^b}}}{t} + \frac{\pi }{2}} \right)} \right| \geqslant {\delta _1} = \mathop {\min }\limits_{1 \leqslant {k_2} \leqslant t - 1} \left| {\sin \left( {\frac{{\pi {k_2}}}{t} + \frac{\pi }{2}} \right)} \right| > 0.\]
\textbf{Лемма 2 доказана.}

Отметим, что условие  (12) выполняется например для оператора $l$  вида (7) (см.[23]).

Следует заметить, что если $\tau  = \frac{a}{\pi }$   является иррациональным числом, то не всегда можно отделить знаменатель выражения  от нуля. Но в некоторых случаях можно найти зависимость  «малости» знаменателя  от номера $k$. Справедлива лемма.

\textbf{Лемма 3.}  Если $\tau  = \frac{a}{\pi } > 0$  является иррациональным алгебраическим числом степени $p \geqslant 2$, то существует число  $L\left( {\tau ,\varepsilon } \right) > 0$  такое , что при всех  $p \in N,{\mkern 1mu} {\kern 1pt} 0 < \varepsilon  < b$ и $k \to  + \infty $   справедлива оценка:
\[\left| {{\Delta _4}\left( k \right)} \right| \geqslant \frac{L}{{{k^{b + b\varepsilon }}}},\eqno(14)\]
при \[{\lambda _k} = {k^{2s}} + O\left( {{k^r}} \right),\,0 < b + b\varepsilon  < 2s - r - b.\]
\[\left| {{\Delta _4}\left( k \right)} \right| \geqslant \frac{L}{{{k^{2s - r - b}}}},\eqno(15)\]
при
\[{\lambda _k} = {k^{2s}} + O\left( {{k^r}} \right),\,0 < 2s - r - b < b + b\varepsilon .\]

\textbf{Доказательство.} Докажем неравенство (14)  ((15) доказывается аналогично). Пусть
\[\left| {{\Delta _4}\left( k \right)} \right| = \left| {\sin \left( {\pi \sqrt[{2n}]{{{\lambda _k}}}\tau  + \frac{\pi }{2}} \right)} \right|,\]
отсюда имеем
\[\left| {{\Delta _4}\left( k \right)} \right| \geqslant \left| {\sin \left( {\pi {k^b}\tau  + \frac{\pi }{2}} \right)} \right|\left| {\cos \left( {aO\left( {{k^{r + b - 2s}}} \right)} \right)} \right| - \left| {\sin \left( {aO\left( {{k^{r + b - 2s}}} \right)} \right)} \right|\left| {sin\left( {\pi {k^b}\tau } \right)} \right| = \]
\[ = \left| {\sin \left( {\pi {k^b}\tau  + \frac{\pi }{2}} \right)} \right|\left| {\left| {\cos \left( {aO\left( {{k^{r + b - 2s}}} \right)} \right)} \right| - \left| {sin\left( {\pi {k^b}\tau } \right)} \right|\left| {\frac{{\sin \left( {aO\left( {{k^{r + b - 2s}}} \right)} \right)}}{{\left| {\sin \left( {\pi {k^b}\tau  + \frac{\pi }{2}} \right)} \right|}}} \right|} \right|.\eqno(16)\]
Используя выпуклость функции $y = \sin x$  на интервале $\left( {0,\frac{\pi }{2}} \right)$ имеем неравенство
\[\left| {\sin x} \right| \geqslant \frac{{2\left| x \right|}}{\pi },{\mkern 1mu} {\kern 1pt} \left| x \right| \leqslant \frac{\pi }{2},\]
далее
\[\left| {\sin \left( {\pi \tau {k^b} + \frac{\pi }{2} - \pi m} \right)} \right| = \left| {\sin \pi {k^b}\left( {\tau  - \frac{{2m - 1}}{{2{k^b}}}} \right)} \right|,\]
где $m \in N$-произвольно. Теперь подберем $m$ так, чтобы выполнялось неравенство
\[\left| {\tau  - \frac{{2m - 1}}{{2{k^b}}}} \right| \leqslant \frac{1}{{2{k^b}}},\]
для этого достаточно положить
\[m = \left[ {\tau {k^b}} \right] + 1,\]
где $\left[ {\tau {k^b}} \right]$  - целая часть  числа $\tau {k^b}$  . Теперь имеем
\[\left| {\sin \pi {k^b}\left( {\tau  - \frac{{2m - 1}}{{2{k^b}}}} \right)} \right| \geqslant \frac{2}{\pi }\left| {\pi {k^b}\left( {\tau  - \frac{{2m - 1}}{{2{k^b}}}} \right)} \right| = 2{k^b}\left| {\tau  - \frac{{2m - 1}}{{2{k^b}}}} \right|.\]
Известно [26], что для любого алгебраического числа $\tau $  степени $p \geqslant 2$  и произвольного $0 < \varepsilon  < 1$ найдется  $\delta \left( {\tau ,\varepsilon } \right) > 0$ такое, что для любой рациональной дроби $\frac{s}{q}$  выполняется неравенство
\[\left| {\tau  - \frac{s}{q}} \right| \geqslant \frac{{\delta \left( {\tau ,\varepsilon } \right)}}{{{q^{2 + \varepsilon }}}}.\]
Используя это  имеем
\[\left| {\sin \pi {k^b}\left( {\tau  - \frac{{2m - 1}}{{2{k^b}}}} \right)} \right| \geqslant 2{k^b}\left| {\tau  - \frac{{2m - 1}}{{2{k^b}}}} \right| \geqslant 2{k^b}\frac{\delta }{{{{\left( {2{k^b}} \right)}^{2 + \varepsilon }}}} = \frac{L}{{{k^{b + b\varepsilon }}}}.\eqno(17)\]
Подставляя (17) в (16) получим требуемый результат. Остальные случаи для выражения ${\Delta _4}\left( k \right)$  рассматриваются аналогичным образом.\textbf{ Лемма 3  доказана}.

Учитывая вышесказанное, получим условия, при которых ряд
\[u\left( {x,y} \right) = \sum\limits_{k = 1}^\infty  {{Y_k}\left( y \right){X_k}\left( x \right)},\eqno(18)\]
является классическим решением поставленной задачи. Формально имеем
\[\left| {D_y^{2n}u\left( {x,y} \right)} \right| \leqslant L\sum\limits_{k = 1}^\infty  {{\lambda _k}\left| {{X_k}} \right|} \frac{{\sum\limits_{r = 0}^{n - 1} {\left\{ {\left| {{\varphi _{rk}}} \right| + \left| {{\psi _{rk}}} \right|} \right\}} }}{{\left| {{\Delta _4}\left( k \right) + {\Delta _5}\left( k \right)} \right|}},\]
теперь, если ${\Delta _1}\left( k \right) \ne 0,{\mkern 1mu} {\kern 1pt} {\mkern 1mu} {\kern 1pt} \frac{a}{\pi } - $  удовлетворяет условиям леммы 2 , то имеем оценку
\[\left| {D_y^{2n}u\left( {x,y} \right)} \right| \leqslant L\sum\limits_{k = 1}^\infty  {{\lambda _k}\left| {{X_k}\left( x \right)} \right|\sum\limits_{r = 0}^{n - 1} {\left\{ {\left| {{\varphi _{rk}}} \right| + \left| {{\psi _{rk}}} \right|} \right\}}},\eqno(19)\]
если ${\Delta _1}\left( k \right) \ne 0$ и выполняется оценка (14), то
\[\left| {D_y^{2n}u\left( {x,y} \right)} \right| \leqslant L\sum\limits_{k = 1}^\infty  {{\lambda _k}{k^{b + b\varepsilon }}\left| {{X_k}} \right|} \sum\limits_{r = 0}^{n - 1} {\left\{ {\left| {{\varphi _{rk}}} \right| + \left| {{\psi _{rk}}} \right|} \right\}}.\eqno(20)\]
если ${\Delta _1}\left( k \right) \ne 0$ и выполняется оценка (15), то
\[\left| {D_y^{2n}u\left( {x,y} \right)} \right| \leqslant L\sum\limits_{k = 1}^\infty  {{\lambda _k}{k^{2s - r - b}}\left| {{X_k}} \right|} \sum\limits_{r = 0}^{n - 1} {\left\{ {\left| {{\varphi _{rk}}} \right| + \left| {{\psi _{rk}}} \right|} \right\}},\eqno(21)\]
Осталось наложить условия на граничные функции.

\textbf{Теорема 2.}  Пусть выполнены следующие условия:

1. ${\Delta _1}\left( k \right) \ne 0,\forall k;$

2. выполнены условия леммы 2 ;

3. ${\varphi _j}\left( x \right),{\psi _j}\left( x \right) \in {C^{4s}}\left[ {0;\pi } \right];$

4. $$\varphi _j^{2m}\left( 0 \right) = \varphi _j^{2m}\left( \pi  \right) = \psi _j^{2m}\left( 0 \right) = \psi _j^{2m}\left( \pi  \right) = 0,$$ $$\overline \varphi  _j^{2m}\left( 0 \right) = \overline \varphi  _j^{2m}\left( \pi  \right) = \overline \psi  _j^{2m}\left( 0 \right) = \overline \psi  _j^{2m}\left( \pi  \right) = 0,$$ $${\overline \varphi  _j}\left( x \right) = l\left( {{\varphi _j}\left( x \right)} \right),\,{\overline \psi  _j}\left( x \right) = l\left( {{\psi _j}\left( x \right)} \right),$$ $$j = 0,1,...n - 1,{\mkern 1mu} {\kern 1pt} {\mkern 1mu} {\kern 1pt} {\mkern 1mu} {\kern 1pt} m = 0,1,...s - 1.$$

Тогда ряд (18) является классическим решением задачи  $D$.

\textbf{Доказательство.} Покажем сходимость каждого слагаемого в оценке (19). Имеем
\[L\sum\limits_{k = 1}^\infty  {{\lambda _k}\left| {{X_k}\left( x \right)} \right|\left| {{\varphi _{0k}}} \right|}  \leqslant L\sqrt {\sum\limits_{k = 1}^\infty  {{{\left( {\frac{{{X_k}\left( x \right)}}{{{\lambda _k}}}} \right)}^2}} } \sqrt {\sum\limits_{k = 1}^\infty  {{{\left( {\lambda _k^2{\varphi _{0k}}} \right)}^2}} }\eqno(22).\]
Первый множитель сходится в силу неравенства Бесселя. Изучим второй множитель, имеем
\[{\varphi _{0k}} = \int\limits_0^l {{\varphi _0}\left( x \right){X_k}\left( x \right)dx}  = \frac{1}{{{\lambda _k}}}\int\limits_0^l {{\varphi _0}\left( x \right)l\left( {{X_k}} \right)dx}  = \frac{1}{{{\lambda _k}}}\int\limits_0^l {l\left( {{\varphi _0}} \right){X_k}dx}  = \]
\[ = \frac{1}{{\lambda _k^2}}\int\limits_0^l {l\left( {{\varphi _0}} \right)l\left( {{X_k}} \right)dx}=\frac{1}{{\lambda _k^2}}\int\limits_0^l {{l^2}\left( {{\varphi _0}} \right){X_k}dx}.\]
Теперь применив неравенство Бесселя получим сходимость второго множителя в соотношении (22). Итак ряд (22) сходится абсолютно и равномерно.
Сходимость остальных слагаемых в (19) показывается аналогичным образом. \textbf{Теорема 2 доказана.}

\textbf{Теорема 3.}  Пусть выполнены следующие условия:

1. ${\Delta _1}\left( k \right) \ne 0,\forall k;$

2. выполнены условия леммы 3 ;

3. $1 \leqslant b \leqslant s - 1;$

4. ${\varphi _j}\left( x \right),{\psi _j}\left( x \right) \in {C^{4s}}\left[ {0;\pi } \right];$

5. $$\varphi _j^{2m}\left( 0 \right) = \varphi _j^{2m}\left( \pi  \right) = \psi _j^{2m}\left( 0 \right) = \psi _j^{2m}\left( \pi  \right) = 0,$$ $$\overline \varphi  _j^{2m}\left( 0 \right) = \overline \varphi  _j^{2m}\left( \pi  \right) = \overline \psi  _j^{2m}\left( 0 \right) = \overline \psi  _j^{2m}\left( \pi  \right) = 0,$$ $${\overline \varphi  _j}\left( x \right) = l\left( {{\varphi _j}\left( x \right)} \right),\,{\overline \psi  _j}\left( x \right) = l\left( {{\psi _j}\left( x \right)} \right),$$ $$j = 0,1,...n - 1,{\mkern 1mu} {\kern 1pt} {\mkern 1mu} {\kern 1pt} {\mkern 1mu} {\kern 1pt} m = 0,1,...s - 1.$$

Тогда ряд (18) является классическим решением задачи  $D$.

\textbf{Доказательство.} Покажем сходимость каждого слагаемого в оценке (20). Имеем
\[L\sum\limits_{k = 1}^\infty  {{\lambda _k}{k^{b + b\varepsilon }}\left| {{X_k}\left( x \right)} \right|} \left| {{\varphi _{0k}}} \right| \leqslant L\sqrt {\sum\limits_{k = 1}^\infty  {{{\left( {\frac{{{X_k}}}{{\sqrt {{\lambda _k}} }}} \right)}^2}} } \sqrt {\sum\limits_{k = 1}^\infty  {\lambda _k^3{k^{2b + 2b\varepsilon }}\varphi _{0k}^2} }  \leqslant \]
\[ \leqslant L\sqrt {\sum\limits_{k = 1}^\infty  {{{\left( {\frac{{{X_k}}}{{\sqrt {{\lambda _k}} }}} \right)}^2}} } \sqrt {\sum\limits_{k = 1}^\infty  {\lambda _k^3{k^{2s - 1 + 2b\varepsilon }}\varphi _{0k}^2} }  \leqslant L\sqrt {\sum\limits_{k = 1}^\infty  {{{\left( {\frac{{{X_k}}}{{\sqrt {{\lambda _k}} }}} \right)}^2}} } \sqrt {\sum\limits_{k = 1}^\infty  {\lambda _k^4\varphi _{0k}^2} } ,\]
Здесь учли асимптотику (9). Первый множитель сходится за счет оценки (8). Сходимость второго множителя показана в доказательстве теоремы 2. Сходимость ряда (21) показывается аналогично. \textbf{Теорема 3 доказана.}

\textbf{Замечание 1.} Если в теореме 3, $b$ будет равняться $s$, то для существования решения поставленной задачи условия гладкости на граничные функции повысятся.

\textbf{Замечание 2.} Если ${\Delta _1}\left( k \right) = 0,$ при некоторых значениях $k = {k_1},{k_2},...{k_p}$, то для
разрешимости задачи $D$ достаточно выполнения условий ортогональности
 $\left( {{\varphi _j}\left( x \right),{X_k}\left( x \right)} \right) = \left( {{\psi _j}\left( x \right),{X_k}\left( x \right)} \right) = 0,\,k = {k_1},...,{k_p},j = 0,..,n - 1.$

\begin{center}
Литература
\end{center}
1. P. G. Bourgin, R. Duffin, “The Dirichlet problem for the vibrating string equation”, Bull.
Amer. Math. Soc., 45:12 (1939), 851–858.\\
2. F. John, “Diriclet problem for a hyperbolic equation”, Amer. J. Math., 63:1 (1941),
141–154.\\
3. С. Л. Соболев, “Пример корректной задачи для уравнения колебания струны с данными на всей границе”, Докл. АН СССР, 109 (1956), 707–709.\\
4. Ю. М. Березанский, “О задаче Дирихле для уравнения колебания струны”, Укр. матем. журн., 12:4 (1960), 363–372.\\
5.П. П. Мосолов, “О задаче Дирихле для уравнений в частных производных”, Изв. вузов. Матем., 1960, № 3, 213–218.\\
6. Б.И.Пташник,  Некорректные граничные задачи для дифференциальных
уравнений с частными производными ,  Наукова Думка, Киев,  1984.\\
7. В. И. Арнольд, “Малые знаменатели. I. Об отображениях окружности на себя”, Изв.
АН СССР. Сер. матем., 25:1 (1961), 21–86.\\
8. Ю. М. Березанский, Разложение по собственным функциям самосопряженных операторов, Наукова Думка, Киев, 1965.\\
9. К. Б. Сабитов, “Задача Дирихле для уравнений с частными производными высоких порядков”, Матем. заметки, 97:2 (2015), 262–276.\\
10. А. В. Бицадзе,  “Некорректность задачи Дирихле для уравнений смешанного типа в смешанных областях “, Докл. АН СССР, 122:2 (1958). 167–170.\\
11. А. М. Нахушев,   “Критерий единственности задачи Дирихле для уравнения смешанного типа в цилиндрической области “, Дифференц. уравнения, 6:1 (1970), 190–191.\\
12. А. П. Солдатов,  “Задача типа Дирихле для уравнения Лаврентьева–Бицадзе. I. Теоремы единственности“, Докл. РАН, 332:6 (1993),696–698.\\
13. А. П. Солдатов,  “Задача типа Дирихле для уравнения Лаврентьева–Бицадзе. II. Теоремы существования“, Докл. РАН, 333:1 (1993),16–18. \\
14. J. R. Cannon,  “A Dirichlet problem for an equation of mixed type with a discontinuous coefficient“, Ann. Mat. Pura Appl., IV. Ser., 61 (1963), 371–377.\\
15. М. М. Хачев , “Задача Дирихле для обобщенного уравнения Лаврентьева–Бицадзе в прямоугольной области “, Дифференц. уравнения, 14:1 (1978), 136–139.\\
16. Р. И.Сохадзе , “О первой краевой задаче для уравнения смешанного типа в прямоугольнике“, Дифференц. уравнения, 19:1 (1983), 127–134.\\
17. К. Б. Сабитов , “Задача Дирихле для уравнений смешанного типа в прямоугольной области“, Докл. РАН, 413:1 (2007), 23–26.\\
18. К. Б. Сабитов , Е. П. Мелишева,  “Задача Дирихле для нагруженного уравнения смешанного типа в прямоугольной области“, Изв. вузов. Матем., 7 (2013), 62–76.\\
19.  К. Б. Сабитов,    В. А. Гущина (Новикова), “Задача  А.А. Дезина для неоднородного уравнения Лаврентьева–Бицадзе“, Изв. вузов. Матем., 3 (2017), 37–50.\\
20. В. И. Жегалов , “Краевая задача для уравнения смешанного типа высшего порядка“, Докл. АН СССР, 136:2 (1961), 274–276.\\
21. М. М. Смирнов , “Краевая задача со смещением для уравнения смешано-составного типа 4-го порядка“, Дифференц. уравнения,11:9 (1975), 1678–1686.\\
22. К. Б. Сабитов,  “О положительности решения неоднородного уравнения смешанного типа высшего порядка“, Изв. вузов. Матем., 3 (2016), 65–71.\\
23. В. А. Садовничий, “О следах обыкновенных дифференциальных операторов высших порядков”, Матем. сб., 72(114):2 (1967), 293–317.\\
24. Ф.Ж. Трикоми, Интегральные уравнения,Москва,изд.Иност.лит.1960.\\
25. М.А. Наймарк, Линейные дифференциальные операторы, Москва,Наука,1969.\\
26. А.А. Бухштаб,  Теория чисел. Москва,Просвещение, 1966.

\end{document}